\documentclass{elsart3-1}

\usepackage{amssymb,amsmath}
\usepackage[english,francais]{babel}

\newtheorem{theorem}{Theorem}[section]

\newtheorem{e-proposition}[theorem]{Proposition}
\newtheorem{corollary}[theorem]{Corollary}
\newtheorem{e-definition}[theorem]{Definition\rm}

\renewcommand{\theequation}{\arabic{equation}}
\newcommand{\R}{\mathbb{R}}
\newcommand{\N}{\mathbb{N}}
\newcommand{\Z}{\mathbb{Z}}
\newcommand{\C}{\mathbb{C}}
\newcommand{\T}{\mathbb{T}}

\newcommand{\dis}{\displaystyle}
\newcommand{\supp}{\textup{supp}}

\def\og{\leavevmode\raise.3ex\hbox{$\scriptscriptstyle\langle\!\langle$~}}
\def\fg{\leavevmode\raise.3ex\hbox{~$\!\scriptscriptstyle\,\rangle\!\rangle$}}

\begin{document}

\begin{frontmatter}

\selectlanguage{english}
\title{Means of algebraic numbers in the unit disk}
\author{Igor E. Pritsker\thanksref{label1}}
\ead{igor@math.okstate.edu}
\ead[url]{http://www.math.okstate.edu/\~{}igor/}
\thanks[label1]{Research was partially supported by NSA, and by the Alexander von Humboldt Foundation.}
\address{Department of Mathematics, Oklahoma State University,
Stillwater, OK 74078, U.S.A.}

\selectlanguage{english}

\begin{abstract}
Schur studied limits of the arithmetic means $s_n$ of zeros for polynomials of degree $n$ with integer coefficients and simple zeros in the closed unit disk. If the leading coefficients are bounded, Schur proved that $\limsup_{n\to\infty} |s_n| \le 1-\sqrt{e}/2.$ We show that $s_n \to 0$, and estimate the rate of convergence by generalizing the Erd\H{o}s-Tur\'an theorem on the distribution of zeros.
{\it To cite this article: I. E. Pritsker, C. R. Acad. Sci. Paris, Ser. I 336 (2003).}

\vskip 0.5\baselineskip

\selectlanguage{francais}
\noindent{\bf R\'esum\'e}
\vskip 0.5\baselineskip
\noindent
{\bf Moyennes de nombres alg\'ebriques dans le disque unit\'e.} Schur a \'etudi\'e les limites des moyennes arithm\'etiques $s_n$ des z\'eros pour les polyn\^omes \`a coefficients entiers de degr\'e $n$ ayant des z\'eros simples dans le disque unit\'e ferm\'e.
Lorsque les coefficients dominants restent born\'es, Schur a d\'emontr\'e que $\limsup_{n\to\infty} |s_n| \le 1-\sqrt{e}/2$. Nous prouvons que $s_n \to 0$. Nous donnons une estimation du taux de
convergence, gr$\rm\hat{a}$ce \`a une g\'{e}n\'{e}ralisation d'un th\'eor\`eme de Erd\H{o}s-Tur\'an sur la distribution des z\'eros. {\it Pour citer cet article~: I. E. Pritsker, C. R. Acad. Sci. Paris, Ser. I 336 (2003).}

\end{abstract}
\end{frontmatter}

\selectlanguage{english}

\section{Schur's problem and equidistribution of zeros} \label{sec1}

Let $\Z_n(D)$ be the set of polynomials of degree $n$ with integer coefficients and all zeros in the closed unit disk $D$. We denote the subset of $\Z_n(D)$ with simple zeros by $\Z_n^1(D)$. Given $M>0$, we write $P_n=a_nz^n + \ldots\in\Z_n^1(D,M)$ if $|a_n|\le M$ and $P_n\in\Z_n^1(D)$ (respectively $P_n\in\Z_n(D,M)$ if $|a_n|\le M$ and $P_n\in\Z_n(D)$). Schur \cite[\S 8]{Sch} studied the limiting behavior of the arithmetic means $s_n$ of zeros for polynomials from $\Z_n^1(D,M)$ as $n\to\infty,$ where $M>0$ is an arbitrary fixed number. He showed that $\limsup_{n\to\infty} |s_n| \le 1-\sqrt{e}/2,$ and remarked that this $\limsup$ is equal to $0$ for {\em monic} polynomials from $\Z_n(D)$ by Kronecker's theorem \cite{Kr}. We prove that $\lim_{n\to\infty} s_n = 0$ for any sequence of polynomials from Schur's class $\Z_n^1(D,M),\ n\in\N.$ This result is obtained as a consequence of the asymptotic equidistribution of zeros near the unit circle. Namely, if $\{\alpha_k\}_{k=1}^n$ are the zeros of $P_n$, we define the counting measure $\tau_n := \frac{1}{n} \sum_{k=1}^n \delta_{\alpha_k}$, where $\delta_{\alpha_k}$ is the unit point mass at $\alpha_k$. Consider the normalized arclength measure $\mu$ on the unit circumference $\T$, with $d\mu(e^{it}):=\frac{1}{2\pi}dt.$ If the $\tau_n$ converge weakly to $\mu$ as $n\to\infty$ ($\tau_n \stackrel{*}{\rightarrow} \mu$) then $\lim_{n\to\infty} s_n = \lim_{n\to\infty} \int z\,d\tau_n(z) = \int z\,d\mu(z) = 0.$ Thus Schur's problem is solved by the following result.

\begin{theorem} \label{thm1.1}
If $P_n(z) = a_nz^n + \ldots\in\Z_n^1(D),\ n\in\N,$ satisfy $\dis\lim_{n\to\infty} |a_n|^{1/n} = 1$, then $\tau_n \stackrel{*}{\rightarrow} \mu$ as $n\to\infty.$
\end{theorem}
Ideas on the equidistribution of zeros date back to Jentzsch and Szeg\H{o}, cf. \cite[Ch. 2]{AB}. They were developed further by Erd\H{o}s and Tur\'an \cite{ET}, and many others; see \cite{AB} for history and additional references. More recently, this topic received renewed attention in number theory, e.g. in the work of Bilu \cite{Bi}. If the leading coefficients of polynomials are bounded, then we can allow even certain multiple zeros. Define the multiplicity of an irreducible factor $Q$ of $P_n$ as the integer $m_n\ge 0$ such that $Q^{m_n}$ divides $P_n$, but $Q^{m_n+1}$ does not divide $P_n$. If a factor $Q$ occurs infinitely often in a sequence $P_n,\ n\in\N,$ then $m_n=o(n)$ means $\lim_{n\to\infty} m_n/n =0.$ If $Q$ is present only in finitely many $P_n$, then $m_n=o(n)$ by definition.

\begin{theorem} \label{thm1.2}
Assume that $P_n\in\Z_n(D,M),\ n\in\N$. If every irreducible factor in the sequence of polynomials $P_n$ has multiplicity $o(n)$, then $\tau_n \stackrel{*}{\rightarrow} \mu$ as $n\to\infty.$
\end{theorem}

\begin{corollary} \label{cor1.3}
If $P_n(z)=a_n\prod_{k=1}^n (z-\alpha_k),\ n\in\N,$ satisfy the assumptions of Theorem \ref{thm1.1} or \ref{thm1.2}, then
\[
\lim_{n\to\infty} \frac{1}{n}\sum_{k=1}^n \alpha_k^m =0, \quad m\in\N.
\]
\end{corollary}
We also show that the norms $\|P_n\|_{\infty}:= \max_{|z|=1} |P_n(z)|$ have at most subexponential growth.
\begin{corollary} \label{cor1.4}
If $P_n,\ n\in\N,$ satisfy the assumptions of Theorem \ref{thm1.1} or Theorem \ref{thm1.2}, then
\[
\lim_{n\to\infty} \|P_n\|_{\infty}^{1/n} = 1.
\]
\end{corollary}
This result is somewhat unexpected, as we have no direct control of the norm or coefficients (except for the leading one). For example, $P_n(z)=(z-1)^n$ has norm $\|P_n\|_{\infty}=2^n$.


We now consider quantitative aspects of the convergence $\tau_n \stackrel{*}{\rightarrow} \mu$. As an application, we obtain estimates of the convergence rate of $s_n$ to $0$ in Schur's problem. A classical result on the distribution of zeros is due to Erd\H{o}s and Tur\'an \cite{ET}. For $P_n(z) = \sum_{k=0}^n a_k z^k$ with $a_k\in\C,$ let $N(\phi_1,\phi_2)$ be the number of zeros in the sector $\{z\in\C:0\le \phi_1 \le \arg(z) \le \phi_2< 2\pi\},$ where $\phi_1 < \phi_2.$ Erd\H{o}s and Tur\'an \cite{ET} proved that
\begin{align} \label{2.1}
\left|\frac{N(\phi_1,\phi_2) }{n} - \frac{\phi_2-\phi_1}{2\pi}\right| \le 16 \sqrt{\frac{1}{n}\log\frac{\|P_n\|_{\infty}}{\sqrt{|a_0 a_n|}}}.
\end{align}
The constant $16$ was improved by Ganelius, and $\|P_n\|_{\infty}$ was replaced by weaker integral norms by Amoroso and Mignotte; see \cite{AB} for more history and references. Our main difficulty in applying \eqref{2.1} to Schur's problem is the absence of an effective estimate for $\|P_n\|_{\infty},\ P_n\in\Z_n^1(D,M)$. We prove a new ``discrepancy" estimate via energy considerations from potential theory. These ideas originated in part in the work of Kleiner, and were developed by Sj\"{o}gren and H\"{u}sing, see \cite[Ch. 5]{AB}. We also use the Mahler measure of a polynomial $P_n(z) = a_n\prod_{k=1}^n (z-\alpha_k)$, defined by $M(P_n) := \exp\left(\frac{1}{2\pi} \int_0^{2\pi} \log |P_n(e^{it})|\,dt\right).$ Note that $M(P_n) = \lim_{p\to 0} \|P_n\|_p$, where $\|P_n\|_p:=\left(\frac{1}{2\pi} \int_0^{2\pi} |P_n(e^{it})|^p\,dt\right)^{1/p},\ p>0$. Jensen's formula readily gives $M(P_n) = |a_n|\prod_{k=1}^n \max(1,|\alpha_k|)$ \cite[p. 3]{Bo}. Hence $M(P_n)=|a_n|\le M$ for any $P_n\in\Z_n(D,M).$

\begin{theorem} \label{thm2.1}
Let $\phi:\C\to\R$ satisfy $|\phi(z)-\phi(t)|\le A|z-t|,\ z,t\in\C,$ and $\supp(\phi)\subset\{z:|z|\le R\}.$ If $P_n(z) = a_n\prod_{k=1}^n (z-\alpha_k)$ is a polynomial with integer coefficients and simple zeros, then
\begin{align} \label{2.2}
\left|\frac{1}{n}\sum_{k=1}^n \phi(\alpha_k) - \int\phi\,d\mu\right| \le A(2R+1) \sqrt{\frac{\log\max(n,M(P_n))}{n}}, \quad n\ge 55.
\end{align}
\end{theorem}
This theorem is related to recent results of Favre and Rivera-Letelier \cite{FR}, obtained in a different setting. Choosing $\phi$ appropriately, we obtain an estimate of the means $s_n$ in Schur's problem.
\begin{corollary} \label{cor2.2}
If $P_n\in\Z_n^1(D,M)$ then
\[
\left|\frac{1}{n}\sum_{k=1}^n \alpha_k\right| \le 8 \sqrt{\frac{\log{n}}{n}},\quad n\ge \max(M,55).
\]
\end{corollary}
We also have an improvement of Corollary \ref{cor1.4} for Schur's class $\Z_n^1(D,M)$.
\begin{corollary} \label{cor2.3}
If $\{P_n\}_{n=1}^{\infty}\in\Z_n^1(D,M)$ then there is some $c>0$ such that $\|P_n\|_{\infty} \le e^{c\sqrt{n}\log{n}}$ as $n\to\infty.$
\end{corollary}

The proof of Theorem \ref{thm2.1} gives a result for arbitrary polynomials with simple zeros, and for any continuous $\phi$ with finite Dirichlet integral $D[\phi]=\iint(\phi_x^2 +\phi_y^2)\,dA$. Moreover, all arguments may be extended to general sets of logarithmic capacity $1$, e.g. to $[-2,2]$. Using the characteristic function $\phi=\chi_E$, we can prove general discrepancy estimates on arbitrary sets, and obtain an Erd\H{o}s-Tur\'an-type theorem. Our results have a number of applications to the problems on integer polynomials considered in \cite{Bo}.

\section{Proofs} \label{sec3}

\noindent{\em Proof of Theorem \ref{thm1.1}.} Observe that the discriminant $\Delta(P_n):=a_n^{2n-2} \prod_{1\le j<k\le n} (\alpha_j-\alpha_k)^2$ is an integer, as a symmetric form in the zeros of $P_n$. Since $P_n$ has simple roots, we have
$\Delta(P_n)\neq 0$ and $|\Delta(P_n)|\ge 1.$ Using weak compactness, we assume that $\tau_n \stackrel{*}{\rightarrow} \tau,$ where $\tau$ is a probability measure on $D$. Let $K_M(x,t) := \min\left(-\log{|x-t|},M\right).$ Since $\tau_n\times\tau_n \stackrel{*}{\rightarrow} \tau\times\tau,$ we obtain for the energy of $\tau$ that
\begin{align*}
I[\tau] &:=-\iint \log|x-t|\,d\tau(x)\,d\tau(t) =
\lim_{M\to\infty} \left( \lim_{n\to\infty} \iint K_M(x,t)\,
d\tau_n(x)\,d\tau_n(t) \right) \\ &= \lim_{M\to\infty} \left(
\lim_{n\to\infty} \left( \frac{1}{n^2} \sum_{j\neq k}
K_M(\alpha_j,\alpha_k) +\frac{M}{n} \right) \right) \le
\lim_{M\to\infty} \left( \liminf_{n\to\infty} \frac{1}{n^2}
\sum_{j\neq k} \log\frac{1}{|\alpha_j-\alpha_k|} \right) \\ &=
\liminf_{n\to\infty} \frac{1}{n^2}
\log\frac{|a_n|^{2n-2}}{\Delta(P_n)} \le \liminf_{n\to\infty} \frac{1}{n^2}\log|a_n|^{2n-2} = 0.
\end{align*}
Thus $I[\tau]\le 0$. But $I[\nu]>0$ for any probability measure $\nu$ on $D$, except for $\mu$ \cite{La}. Hence $\tau=\mu.$

\noindent{\em Proof of Theorem \ref{thm1.2}.} Let $\phi\in C(\C).$ Note that for any $\epsilon>0$ there are finitely many irreducible factors $Q$ in the sequence $P_n$ such that $|\int\phi\,d\tau(Q) - \int\phi\,d\mu| \ge \epsilon,$ where $\tau(Q)$ is the zero counting measure for $Q$. Indeed, if we have an infinite sequence of such $Q_m$, then $\deg(Q_m)\to\infty$, as there are only finitely many $Q_m\in\Z_n(D,M)$ of bounded degree. Hence $\int\phi\,d\tau(Q_m) \to \int\phi\,d\mu$ by Theorem \ref{thm1.1}. Let the number of such exceptional factors $Q_m$ be $N$. Then we have $|n\int\phi\,d\tau_n - n\int\phi\,d\mu| \le N o(n) \max_D |\phi - \int\phi\,d\mu| + (n-N)\epsilon,\ n\in\N.$ Hence $\limsup_{n\to\infty} |\int\phi\,d\tau_n - \int\phi\,d\mu| \le \epsilon,$ and $\lim_{n\to\infty} \int\phi\,d\tau_n = \int\phi\,d\mu$ after letting $\epsilon\to 0.$

\noindent{\em Proof of Corollary \ref{cor1.3}.} Let $\phi(z)=z^m$ and write $\lim_{n\to\infty} \int z^m\,d\tau_n(z) = \int z^m\,d\mu(z) = 0.$

\noindent{\em Proof of Corollary \ref{cor1.4}.} Let $\|P_n\|_{\infty}=|P_n(z_n)|,\ z_n\in D,$ and assume $\lim_{n\to\infty} z_n = z_0\in D$ by compactness. Then $\|P_n\|_{\infty} = \exp\left(\log\left|P_n(z_n)\right|\right) = |a_n| \exp\left(n \int\log|z_n-t|\,d\tau_n(t)\right)$. Since $\tau_n \stackrel{*}{\rightarrow} \mu$, Theorem I.6.8 of \cite{ST} gives $\limsup_{n\to\infty} \|P_n\|_{\infty}^{1/n} \le \exp\left(\int\log|z_0-t|\,d\mu(t)\right) = 1$ \cite[p. 22]{ST}. But $\|P_n\|_{\infty} \ge |a_n| \ge 1$, see \cite[p. 16]{AB}.

\noindent{\em Proof of Theorem \ref{thm2.1}.} Given $r>0$, define the measures $\nu_k^r$ with $d\nu_k^r(\alpha_k + re^{it}) = dt/(2\pi),\ t\in[0,2\pi).$ Let $\tau_n^r:=\frac{1}{n}\sum_{k=1}^n \nu_k^r$, and estimate $\left|\int\phi\,d\tau_n - \int\phi\,d\tau_n^r\right| \le \frac{1}{n}\sum_{k=1}^n \frac{1}{2\pi}\int_0^{2\pi} \left|\phi(\alpha_k) - \phi(\alpha_k + re^{it})\right|\,dt \le \omega_{\phi}(r),$ where $\omega_{\phi}(r):=\sup_{|z-\zeta|\le r} |\phi(z)-\phi(\zeta)|$ is the modulus of continuity of $\phi$.

Let $p_{\nu}(z):=-\int\log|z-t|d\nu(t)$ be the potential of a measure $\nu.$ A direct evaluation gives that $p_{\nu_k^r}(z)=-\log\max(r,|z-\alpha_k|)$ and $p_{\mu}(z) = -\log\max(1,|z|)$ \cite[p. 22]{ST}. Consider $\sigma:=\tau_n^r-\mu,\ \sigma(\C)=0.$ One computes (or see \cite[p. 92]{ST}) that
$d\sigma=-\frac{1}{2\pi}\left(\partial p_{\sigma}/\partial n_+ + \partial p_{\sigma}/\partial n_-\right) ds,$
where $ds$ is the arclength on $\supp(\sigma)=\{|z|=1\}\cup \left( \cup_{k=1}^n \{|z-\alpha_k|=r\}\right)$, and $n_{\pm}$ are the inner and the outer normals. We now use Green's identity $\iint_G u \Delta v\,dA =  \int_{\partial G} u\,\frac{\partial v}{\partial n}\,ds - \iint_G \nabla u \cdot \nabla v\,dA$ with $u=\phi$ and $v=p_{\sigma}$ in each component $G$ of $\{|z|<R\}\setminus\supp(\sigma).$ Since $\Delta p_{\sigma}=0$ in $G$, adding the identities for all $G$, we obtain that
\[
\left|\int\phi\,d\sigma\right| = \frac{1}{2\pi} \left| \iint_{|z|\le R} \nabla \phi \cdot \nabla p_{\sigma} \,dA \right| \le \frac{1}{2\pi} \sqrt{D[\phi]}\,\sqrt{D[p_{\sigma}]},
\]
where $D[\phi]=\iint(\phi_x^2 +\phi_y^2)\,dA$ is the Dirichlet integral of $\phi.$ It is known that $D[p_{\sigma}]=2\pi I[\sigma]$ \cite[Thm 1.20]{La}, where $I[\sigma]=-\iint \log|z-t|\,d\sigma(z)\,d\sigma(t) = \int p_{\sigma}\,d\sigma$. Since $p_{\mu}(z) = -\log\max(1,|z|)$, we observe that $\int p_{\mu}\,d\mu = 0,$ so that
$I[\sigma]=\int p_{\tau_n^r}\,d\tau_n^r - 2\int p_{\mu}\,d\tau_n^r.$ Further, $-\int p_{\mu}\,d\tau_n^r = \int\log\max(1,|z|)\, d\tau_n^r(z) \le \left(\sum_{|\alpha_k|\le 1+r} \log(1+2r) + \sum_{|\alpha_k|>1+r} \log|\alpha_k|\right)/n \le \log(1+2r) +\frac{1}{n}\log M(P_n) - \frac{1}{n}\log |a_n|.$  We also have that
$\int p_{\tau_n^r}\,d\tau_n^r \le \left(-\sum_{j\neq k} \log|\alpha_j-\alpha_k| - n\log{r}\right)/n^2.$ We next combine the energy estimates to obtain
\[
I[\sigma] \le \frac{2}{n}\log M(P_n) - \frac{1}{n^2}\log\left|a_n^2 \Delta(P_n)\right| - \frac{1}{n}\log{r} + 4r.
\]
Collecting all estimates, we proceed with $\left|\int\phi\,d\tau_n - \int\phi\,d\mu\right| \le \left|\int\phi\,d\tau_n - \int\phi\,d\tau_n^r\right| + \left|\int\phi\,d\tau_n^r - \int\phi\,d\mu\right| \le \omega_{\phi}(r) + \sqrt{D[\phi]}\,\sqrt{D[p_{\sigma}]}/(2\pi) = \omega_{\phi}(r) + \sqrt{D[\phi]}\,\sqrt{I[\sigma]/(2\pi)}$. Thus we arrive at the main inequality:
\begin{align} \label{2.3}
\left|\int\phi\,d\tau_n - \int\phi\,d\mu\right| \le \omega_{\phi}(r) + \sqrt{\frac{D[\phi]}{2\pi}}\,\left(\frac{2}{n}\log M(P_n) - \frac{1}{n^2}\log\left|a_n^2 \Delta(P_n)\right| - \frac{1}{n}\log{r} + 4r\right)^{1/2}.
\end{align}
Note that $D[\phi]\le 2\pi R^2 A^2,$ as $|\phi_x|\le A$ and $|\phi_y|\le A$ a.e. in $\C.$ Also, $\omega_{\phi}(r)\le Ar.$ Since $|\Delta(P_n)|\ge 1$ and $|a_n|\ge 1$, we have $|a_n^2\Delta(P_n)|\ge 1$. Hence \eqref{2.2} follows from \eqref{2.3} by letting $r=1/\max(n,M(P_n)).$

\noindent{\em Proof of Corollary \ref{cor2.2}.} Since $P_n$ has real coefficients, we have that $s_n = \int z\,d\tau_n(z) = \int \Re(z)\,d\tau_n(z).$ We let $\phi(z)=\Re(z),\ |z|\le 1;\ \phi(z)=\Re(z)(1-\log|z|),\ 1\le |z|\le e;$ and $\phi(z)=0,\ |z|\ge e.$ An elementary computation shows that $|\phi_x(z)|\le 1$ and $|\phi_y(z)|\le 1/2$ for all $z=x+iy\in\C.$ The Mean Value Theorem gives $|\phi(z)-\phi(t)|\le |z-t| \max_{\C} \sqrt{\phi_x^2+\phi_y^2}.$ Hence we can use Theorem \ref{thm2.1} with $A=\sqrt{5}/2$ and $R=e.$

\noindent{\em Proof of Corollary \ref{cor2.3}.} Note that $\log|P_n(z)| = n \int\log|z-w|\,d\tau_n(w).$ For $|z|=1+1/n,$ we let $\phi(w)=\log|z-w|,\ |w|\le 1;\ \phi(w)=(1-\log|w|)\log|1-\bar{z}w|,\ 1\le |w|\le e;$ and $\phi(z)=0,\ |w|\ge e.$ Then $|\phi_x(w)|=O(|z-w|^{-1}),\ |w|\le 1; |\phi_x(w)|=O(|1-\bar{z}w|^{-1}),\ 1\le |w|\le e;$ and the same estimates hold for $|\phi_y|.$ Hence $D[\phi] = O\left(\iint_{|w|\le 1} |z-w|^{-2} dA(w)\right) =  O\left(\int_{1/n}^1 r^{-1}dr \right) = O(\log{n}),$ and $\omega_\phi(r) \le r \max_{\C} \sqrt{\phi_x^2+\phi_y^2} = r O(n).$ Let $r=1/n^2$ and use \eqref{2.3} to obtain $|\log|P_n(z)| - n\log|z|| = O(\sqrt{n}\log{n}).$


\end{document}